\documentclass[12pt,reqno]{article}
\usepackage{amsmath, amssymb}
\newtheorem{theorem}{Theorem}%[section]
\newtheorem{definition}[theorem]{Definition}
\newtheorem{proposition}[theorem]{Proposition}
\newtheorem{corollary}[theorem]{Corollary}
\newtheorem{lemma}[theorem]{Lemma}

\def \proof {\noindent {\bf Proof.}\ \ }
\def \remark {\noindent {\bf Remark.}\ \ }

\def \endproof
{{\mbox{}\nolinebreak\hfill\rule{2mm}{2mm}\par\medbreak}}

\def \R {\mathbb{R}}

\def \Z {\mathbb{Z}}
\def \E {\mathbb{E}}
\def \P {\mathbb{P}}

\def \NN {\mathcal{N}}

\def \b {\beta}

\def \e {\varepsilon}
\def \eps {\varepsilon}
\def \d {\delta}

\def \l {\lambda}

\def \s {\sigma}

\def \w {\omega}

\def \< {\langle}
\def \> {\rangle}

\def \conv {{\rm conv}}

\def \vc {{\rm vc}}

\begin{document}
\title {Entropy and the Combinatorial Dimension}
\author {S. Mendelson\footnote{Research School of
Information Sciences and Engineering, The Australian National
University, Canberra, ACT 0200, Australia, e-mail:
shahar.mendelson@anu.edu.au} \and
   R. Vershynin\footnote{
   Department of Mathematical Sciences,
   University of Alberta,
   Edmonton, Alberta T6G 2G1, Canada,
   e-mail: vershynin@yahoo.com}}
\date{}

\maketitle
%\tableofcontents
%\newpage

\begin{abstract}
We solve Talagrand's entropy problem: the $L_2$-covering numbers
of every uniformly bounded class of functions are exponential in
its shattering dimension. This extends Dudley's theorem on classes
of $\{0,1\}$-valued functions, for which the shattering dimension
is the Vapnik-Chervonenkis dimension.

In convex geometry, the solution means that the entropy of a
convex body $K$ is controlled by the maximal dimension of a cube
of a fixed side contained in the coordinate projections of $K$.
This has a number of consequences, including the optimal Elton's
Theorem and estimates on the uniform central limit theorem in the
real valued case.
\end{abstract}

%MSC Classification: 47A30, 42A60, 94A12

\section{Introduction}
%----------------------------------------------------------------------

The fact that the covering numbers of a set are exponential in its
linear algebraic dimension is fundamental and simple. Let $A$ be a
class of functions bounded by $1$, defined on a set $\Omega$. If
$A$ is a finite dimensional class then for every probability
measure on $\mu$ on $\Omega$,
\begin{equation}                                \label{volumetric}
  N(A, t, L_2(\mu))  \le \Big( \frac{3}{t} \Big)^{\dim(A)},
  \ \ \ \ 0 < t < 1,
\end{equation}
where $\dim(A)$ is the linear algebraic dimension of $A$ and the
left-hand side of \eqref{volumetric} is the covering number of
$A$, the minimal number of functions needed to approximate any
function in $A$ within an error $t$ in the $L_2(\mu)$-norm. This
inequality follows by a simple volumetric argument (see e.g.
\cite{Pi} Lemma 4.10) and is, in a sense, optimal: the dependence
both on $t$ and on the dimension is sharp (except, perhaps, for
the constant $3$).

The linear algebraic dimension of $A$ is often too large for
\eqref{volumetric} to be useful, as it does not capture the
``size'' of $A$ in different directions but only determines in how
many directions $A$ does not vanish. The aim of this paper is to
replace the linear algebraic dimension by a combinatorial
dimension originated from the classical works of Vapnik and
Chervonenkis \cite{VC 71}, \cite{VC 81}.

We say that a subset $\s$ of $\Omega$ is $t$-shattered by a class
$A$ if there exists a level function $h$ on $\s$ such that, given
any subset $\s'$ of $\s$, one can find a function $f \in A$ with
$f(x)  \le  h(x)$ if $x \in \s'$ and $f(x)  \ge  h(x) + t$ if $x
\in \s \setminus \s'$. The {\em shattering dimension} of $A$,
denoted by $\vc(A,t)$ after Vapnik and Chervonenkis, is the
maximal cardinality of a set $t$-shattered by $A$. Clearly, the
shattering dimension does not exceed the linear algebraic
dimension, and is often much smaller. Our main result states that
the linear algebraic dimension in \eqref{volumetric} can be
essentially replaced by the shattering dimension.

\begin{theorem}                                       \label{main}
  Let $A$ be a class of functions bounded by $1$,
  defined on a set $\Omega$.
  Then for every probability measure $\mu$ on $\Omega$,
  \begin{equation}                          \label{combinatorial}
    N(A, t, L_2(\mu))
    \le  \Big( \frac{2}{t} \Big)^{K \cdot \vc(A,\, c t)},
     \ \ \ \ 0 < t < 1,
  \end{equation}
  where $K$ and $c$ are positive absolute constants.
\end{theorem}
There also exists a (simple) reverse inequality complementing
\eqref{combinatorial}: for some measure $\mu$, one has $N(A, t,
L_2(\mu))  \ge  2^{K \cdot \vc(A,\, c t)}$, where $K$ and $c$ are
some absolute constants, see e.g. \cite{T 02}.

The origins of Theorem \ref{main} are rooted in the work of Vapnik
and Chervonenkis, who first understood that entropy estimates are
essential in determining whether a class of functions obeys the
uniform law of large numbers. The subsequent fundamental works of
Koltchinskii \cite{K} and Gin\`{e} and Zinn \cite{GZ} enhanced the
link between entropy estimates and uniform limit theorems (see
also \cite{T 96}).

In 1978, R.~Dudley proved Theorem \ref{main} for classes of
$\{0,1\}$-valued functions (\cite{Du}, see \cite{LT} 14.3). This
yielded that a $\{0,1\}$-class obeys the uniform law of large
numbers (and even the uniform Central Limit Theorem) if and only
if its shattering dimension is finite for $0 < t < 1$. The main
difficulty in proving such limit theorems for general classes has
been the absence of a uniform entropy estimate of the nature of
Theorem \ref{main} (\cite{T 88}, \cite{T 92}, \cite{T 96},
\cite{ABCH}, \cite{BL}, \cite{T 02}). However, proving Dudley's
result for general classes is considerably more difficult due to
the lack of the obvious property of the $\{0,1\}$-valued classes,
namely that if a set $\s$ is $t$-shattered for some $0 < t < 1$
then it is automatically $1$-shattered.

In 1992, M.~Talagrand proved a weaker version of Theorem
\ref{main}: under some mild regularity assumptions, $\log N(A, t,
L_2(\mu)) \le  K \cdot \vc(A, \,ct) \log^M (\frac{2}{t})$, where
$K$, $c$ and $M$ are some absolute constants (\cite{T 92}, \cite{T
02}). Theorem \ref{main} is Talagrand's inequality with the best
possible exponent $M = 1$ (and without regularity assumptions).

Talagrand's inequality was motivated not only by limit theorems in
probability, but to a great extent by applications to convex
geometry. A subset $B$ of $\R^n$ can be viewed as a class of real
valued functions on $\{1, \ldots, n\}$. If $B$ is convex and, for
simplicity, symmetric, then its shattering dimension $\vc(B,t)$ is
the maximal cardinality of a subset $\s$ of $\{1, \ldots, n\}$
such that $P_\s (B)  \supset [-\frac{t}{2}, \frac{t}{2}]^\s$,
where $P_\s$ denotes the orthogonal projection in $\R^n$ onto
$\R^\s$. In the general, non-symmetric, case we allow translations
of the cube $[-\frac{t}{2}, \frac{t}{2}]^\s$ by a vector in
$\R^\s$.

The following entropy bound for convex bodies is then an immediate
consequence of Theorem \ref{main}. Recall that $N(B,D)$ is the
covering number of $B$ by a set $D$ in $\R^n$, the minimal number
of translates of $D$ needed to cover $B$.

\begin{corollary}                 \label{main corollary}
  There exist positive absolute constants $K$ and $c$ such that
  the following holds.
  Let $B$ be a convex body contained in $[0,1]^n$,
  and $D_n$ be the unit Euclidean ball in $\R^n$.
  Then for $0 < t < 1$
  $$
  N(B, t \sqrt{n} D_n)
  \le  \Big( \frac{2}{t} \Big)^{K d},
  $$
  where $d$ is the maximal cardinality of a subset $\s$
  of $\{1, \ldots, n\}$ such that
  $$
  P_\s (B)  \supseteq  h + [0, ct]^\s
  \ \ \ \text{for some vector $h$ in $\R^n$.}
  $$
\end{corollary}

As M.~Talagrand notices in \cite{T 02}, Theorem \ref{main} is a
``concentration of pathology'' phenomenon. Assume one knows that a
covering number of the class $A$ is large. All this means is that
$A$ contains many well separated functions, but it tells nothing
about the structure these functions form. The conclusion of
\eqref{combinatorial} is that $A$ must shatter a large set $\s$,
which detects a very accurate pattern: one can find functions in
$A$ oscillating on $\s$ in all possible $2^{|\s|}$ ways around
fixed levels. The ``largeness'' of $A$, {\em a priori} diffused,
is {\em a fortiori} concentrated on the set $\s$.

The same phenomenon is seen in Corollary \ref{main corollary}:
given a convex body $B$ with large entropy, one can find an entire
cube in a coordinate projection of $B$, the cube that certainly
{\em witnesses} the entropy's largeness.

When dualized, Corollary \ref{main corollary} solves the problem
of finding the best asymptotics in Elton's Theorem. Let $x_1,
\ldots, x_n$ be vectors in the unit ball of a Banach space, and
$\e_1, \ldots, \e_n$ be Rademacher random variables (independent
Bernoulli random variables taking values $1$ and $-1$ with
probability $1/2$). By the triangle inequality, the expectation
$\E \| \sum_{i=1}^n \e_i x_i\|$ is at most $n$, and assume that
$\E \| \sum_{i=1}^n \e_i x_i \|  \ge  \d n$ for some number $\d >
0$.

In 1983, J.~Elton \cite{E} proved an important result that there
exists a subset $\s$ of $\{1, \ldots, n\}$ of size proportional to
$n$ such that the set of vectors $(x_i)_{i \in \s}$ is equivalent
to the $\ell_1$ unit-vector basis. Specifically, there exist
numbers $s, t > 0$, depending only on $\d$, such that
\begin{equation}                                \label{st}
|\s| \ge s^2 n \ \ \ \text{and} \ \ \ \Big\| \sum_{i \in \s} a_i
x_i \Big\|
  \ge  t \sum_{i \in \s} |a_i|
\ \ \text{for all real numbers $(a_i)$}.
\end{equation}
Several steps have been made towards finding the best possible $s$
and $t$ in Elton's Theorem. A trivial upper bound is $s, t \le \d$
which follows from the example of identical vectors and by
shrinking the usual $\ell_1$ unit-vector basis. As for the lower
bounds, J.~Elton proved \eqref{st} with $s \sim \d / \log(1/\d)$
and $t \sim \d^3$. A.~Pajor \cite{Pa} removed the logarithmic
factor from $s$. M.~Talagrand \cite{T 92}, using his inequality
discussed above, improved $t$ to $\d / \log^M (1/\d)$. In the
present paper, we use Corollary \ref{main corollary} to solve this
problem by proving the optimal asymptotics: $s, t \sim \d$.

\begin{theorem}                             \label{intro elton}
   Let $x_1, \ldots, x_n$ be vectors in the unit
   ball of a Banach space, satisfying
   $$
   \E \Big\| \sum_{i=1}^n \e_i x_i \Big\|  \ge  \d n
   \ \ \ \text{for some number $\d > 0$}.
   $$
   Then there exists a subset $\s \subset \{ 1, \ldots, n \}$
   of cardinality $|\s|  \ge  c \d^2 n$ such that
   $$
   \Big\| \sum_{i \in \s} a_i x_i \Big\|
   \ge   c \d \sum_{i \in \s} |a_i|
   \ \ \text{for all real numbers $(a_i)$},
   $$
   where $c$ is a positive absolute constant.
\end{theorem}
Furthermore, there is an interplay between the size of $\s$ and
the isomorphism constant -- they can not attain their worst
possible values together. Namely, we prove that $s$ and $t$ in
\eqref{st} satisfy in addition to $s,t \gtrsim \d$ also the lower
bound $s \cdot t \log^{1.6}(2/t) \gtrsim \d$, which, as an easy
example shows, is optimal for all $\d$ within the logarithmic
factor. The power 1.6 can be replaced by any number greater than
1.5. This estimate improves one of the main results of the paper
\cite{T 92} where this phenomenon in Elton's Theorem was
discovered and proved with a constant (unspecified) power of
logarithm.

The paper is organized as follows. In the remaining part of the
introduction we sketch the proof of Theorem \ref{main}; the
complete proof will occupy Section \ref{s:proof}. Section
\ref{s:convexity} is devoted to applications to Elton's Theorem
and to empirical processes.

Here is a sketch of the proof of Theorem \ref{main}. Starting with
a set $A$ which is separated with respect to the $L_2(\mu)$-norm,
it is possible find a coordinate $\omega \in \Omega$ (selected
randomly) on which $A$ is diffused, i.e. the values $\{f(\omega),
\; f \in A\}$ are spread in the interval $[-1,1]$. Then there
exist two nontrivial subsets $A_1$ and $A_2$ of $A$ with their set
of values $\{f(\omega), \; f \in A_1\}$ and $\{f(\omega), \; f \in
A_2\}$ well separated from each other on the line. Continuing this
process of separation for $A_1$ and $A_2$, etc., one can construct
a dyadic tree of subsets of $A$, called a separating tree, with at
least $|A|^{1/2}$ leaves. The ``largeness'' of the class $A$ is
thus captured by its separating tree.

The next step evoked from a beautiful idea in \cite{ABCH}. First,
there is no loss of generality in discretizing the class: one can
assume that $\Omega$ is finite (say $|\Omega|=n$) and that the
functions in $A$ take values in $\frac{t}{6}\Z \cap [-1,1]$. Then,
instead of producing a large set $\s$ shattered by $A$ with a
certain level function $h$, one can count the number of different
pairs $(\s,h)$ for which $\s$ is shattered by $A$ with the level
function $h$. If this number exceeds $\sum_{k=0}^d
\binom{n}{k}(\frac{12}{t})^k$ then there must exist a set $\s$ of
size $|\s|>d$ shattered by $A$ (because there are $\binom{n}{k}$
possible sets $\s$ of cardinality $k$, and for such a set there
are at most $(\frac{12}{t})^k$ possible level functions).

The only thing remaining is to bound below the number of pairs
$(\s,h)$ for which $\s$ is shattered by $A$ with a level function
$h$. One can show that this number is bounded below by the number
of the leaves in the separating tree of $A$, which is $|A|^{1/2}$.
This implies that $|A|^{1/2} \leq\sum_{k=0}^d
\binom{n}{k}(\frac{12}{t})^k \sim (\frac{n}{td})^d$, where
$d=\vc(A,ct)$. The ratio $\frac{n}{d}$ can be eliminated from this
estimate by a probabilistic extraction principle which reduces the
cardinality of $\Omega$.

\newpage

\noindent ACKNOWLEDGEMENTS

\noindent The first author was supported by an Australian Research
Council Discovery grant. The second author thanks Nicole
Tomczak-Jaegermann for her constant support. He also acknowledges
a support from the Pacific Institute of Mathematical Sciences, and
thanks the Department of Mathematical Sciences of the University
of Alberta for its hospitality. Finally, we would like to thank
the referee for his valuable comments and suggestions.

\section{The Proof of Theorem \ref{main}}
\label{s:proof}
%----------------------------------------------------------------------

For $t > 0$, a pair of functions $f$ and $g$ on $\Omega$ is {\em
$t$-separated in $L_2(\mu)$} if $\|f - g\|_{L_2(\mu)} > t$. A set
of functions is called $t$-separated if every pair of distinct
points in the set is $t$-separated. Let $N_{\rm sep} (A, t,
L_2(\mu))$ denote the maximal cardinality of a $t$-separated
subset of $A$. It is standard and easily seen that
\begin{equation*}
  N (A, t, L_2(\mu))
  \le N_{\rm sep} (A, t, L_2(\mu))
  \le N (A, \frac{t}{2}, L_2(\mu)).
\end{equation*}
This inequality shows that in the proof of Theorem \ref{main} we
may assume that $A$ is $t$-separated in the $L_2(\mu)$ norm, and
replace its covering number by its cardinality.

We will need two probabilistic results, the first of which is
straightforward.

\begin{lemma}                                  \label{l:xx'}
  Let $X$ be a random variable
  and $X'$ be an independent copy of $X$.
  Then
  $$
  \E |X - X'|^2 = 2 \E |X - \E X|^2 = 2 \inf_a \E |X - a|^2.
  $$
\end{lemma}

The next lemma is a small deviation principle. Denote by $\s(X)^2
= \E |X - \E X|^2$ the variance of the random variable $X$.

\begin{lemma}                               \label{l:small deviation}
  Let $X$ be a random variable with nonzero variance.
  Then there exist numbers $a \in \R$ and
  $0 < \b \le \frac{1}{2}$, so that letting
  \begin{align*}
  p_1 &= \P \{ X > a + {\textstyle \frac{1}{6}} \s(X) \}
  \ \ \ \text{and} \\
  p_2 &= \P \{ X < a - {\textstyle \frac{1}{6}} \s(X) \},
  \end{align*}
  one has either $p_1 \ge 1-\b$ and $p_2 \ge \frac{\b}{2}$,
  or $p_2 \ge 1-\b$ and $p_1 \ge \frac{\b}{2}$.
\end{lemma}

\proof Recall that a median of $X$ is a number $M_X$ such that $\P
\{ X \ge M_X \} \ge 1/2$ and $\P \{ X \le M_X \} \ge 1/2$; without
loss of generality we may assume that $M_X = 0$. Therefore $\P \{
X > 0 \} = 1 - \P \{ X \le 0 \} \le  1/2$ and similarly $\P \{ X <
0 \} \le 1/2$.

By Lemma \ref{l:xx'},
\begin{align}                                        \label{integrals}
\s(X)^2
  &\le \E |X|^2
   =   \int_0^\infty \P \{ |X| > \l \} \; d\l^2 \nonumber \\
  &=   \int_0^\infty \P \{ X > \l \} \; d\l^2
     + \int_0^\infty \P \{ X < -\l \} \; d\l^2
\end{align}
where $d\lambda^2 = 2\lambda \;d\lambda$.

Assume that the conclusion of the lemma fails, and let $c$ be any
number satisfying $\frac{1}{3} < c < \frac{1}{\sqrt{8}}$. Divide
$\R_+$ into intervals $I_k$ of length $c \s(X)$ by setting
$$
I_k  =  \Big( c \s(X) k, \; c \s(X) (k + 1) \Big], \ \ \ k = 0, 1,
2, \ldots
$$
and let $\b_0, \b_1, \b_2, \ldots$ be the non-negative numbers
defined by
$$
\P \{ X > 0 \}    = \b_0 \le 1/2, \ \ \ \ \ P \{X \in I_k \} =
\b_k - \b_{k+1}, \ \ \ k = 0, 1, 2, \ldots
$$
We claim that
\begin{equation}                                           \label{bk}
  \text{for all $k \ge 0$,} \ \ \
  \b_{k+1}  \le  \frac{1}{2} \b_k.
\end{equation}
Indeed, assume that $\b_{k+1} > \frac{1}{2} \b_k$ for some $k$ and
consider the intervals $J_1 = \left(-\infty, c \s(X) k \right]$
and $J_2 = \left(c \s(X) (k+1), \infty \right)$. Then $J_1 =
(-\infty, 0] \cup ( \bigcup_{0 \le l \le k-1} I_l)$, so
$$
\P \{ X \in J_1 \} = (1 - \b_0) + \sum_{0 \le l \le k-1} (\b_l -
\b_{l+1}) = 1- \b_k.
$$
Similarly, $J_2 = \bigcup_{l \ge k+1} I_l$ and thus
$$
\P \{ X \in J_2 \} = \sum_{l \ge k+1} (\b_l - \b_{l+1}) = \b_{k+1}
> \frac{1}{2} \b_k.
$$
Moreover, since the sequence $(\b_k)$ is non-increasing by its
definition, then $\b_k \ge \b_{k+1} > \frac{1}{2} \b_k \ge 0$ and
$\b_k \le \b_0 \le \frac{1}{2}$. Then the conclusion of the lemma
would hold with $a$ being the middle point between the intervals
$J_1$ and $J_2$ and with $\b = \b_k$, which contradicts the
assumption that the conclusion of the lemma fails. This proves
\eqref{bk}.

Now, one can apply \eqref{bk} to estimate the first integral in
\eqref{integrals}. Note that whenever $\l \in I_k$,
$$
\P \{ X > \l \} \le \P \{ X > c \s(X) k \} =  \P \Big( \bigcup_{l
\ge k} I_l \Big) =  \b_k.
$$
Then
\begin{align}
\int_0^\infty \P \{ X > \l \} \; d\l^2
  &\le  \sum_{k \ge 0} \int_{I_k} \b_k \cdot 2 \l \;d\l \nonumber\\
  &\le  \sum_{k \ge 0} \b_k \cdot
        2 c \s(X) (k+1) \;{\rm length}(I_k).          \label{series}
\end{align}
Applying \eqref{bk} inductively, it is evident that $\b_k \le
(\frac{1}{2})^k \b_0 \le  \frac{1}{2^{k+1}}$, and since ${\rm
length }(I_k) = c \s(X)$, \eqref{series}  is bounded by
$$
2 c^2 \s(X)^2 \sum_{k \ge 0} \frac{k+1}{2^{k+1}} = 4 c^2 \s(X)^2 <
\frac{1}{2} \s(X)^2.
$$

By an identical argument one can show that the second integral in
\eqref{integrals} is also bounded by $\frac{1}{2} \s(X)^2$.
Therefore
$$
\s(X)^2  < \frac{1}{2} \s(X)^2 + \frac{1}{2} \s(X)^2 = \s(X)^2,
$$
and this contradiction completes the proof.
\endproof

\subsection*{Constructing a separating tree}
%......................................................................

Let $A$ be a finite class of functions on a probability space
$(\Omega, \mu)$, which is $t$-separated in $L_2 (\mu)$. Throughout
the proof we will assume that $|A| > 1$. One can think of the
class $A$ itself as a (finite) probability space with the uniform
measure on it, that is, each element $x$ in $A$ is assigned
probability $\frac{1}{|A|}$.

\begin{lemma}                       \label{separationlemma}
  Let $A$ be a $t$-separated subset of $L_2(\mu)$. Then, there
  exist a coordinate $i$ in $\Omega$
  and numbers $a \in \R$ and $0<\beta \leq 1/2$, so that setting
  \begin{align*}
    N_1 &= | \{ x \in A : \; x(i) > a + {\textstyle \frac{1}{12}} t \}
|
    \ \ \ \text{and} \\
    N_2 &= | \{ x \in A : \; x(i) < a - {\textstyle \frac{1}{12}} t \}
|,
  \end{align*}
  one has either $N_1 \ge  (1 - \b) |A|$
  and $N_2 \ge  \frac{\b}{2} |A|$, or vice versa.
\end{lemma}

\begin{proof}
Let $x, x'$ be random points in $A$ selected independently
according to the uniform (counting) measure on $A$. By Lemma
\ref{l:xx'},
\begin{align}                                         \label{exys}
\E \|x - x'\|_{L_2(\mu)}^2
  &=   \E \int_\Omega |x(i) - x'(i)|^2 \; d\mu(i)
   =   \int_\Omega \E |x(i) - x'(i)|^2 \; d\mu(i)  \nonumber \\
  &= 2 \int_\Omega \E |x(i) - \E x(i)|^2 \; d\mu(i) \\ \nonumber
  &= 2 \int_\Omega \s(x(i))^2 \; d\mu(i)
\end{align}
where $\s(x(i))^2$ is the variance of the random variable $x(i)$
with respect to the uniform measure on $A$.

On the other hand, with probability $1 - \frac{1}{|A|}$ we have $x
\ne x'$ and, whenever this event occurs, the separation assumption
on $A$ implies that $\|x - x'\|_{L_2(\mu)} \ge t$. Therefore
$$
\E \|x - x'\|^2_{L_2(\mu)} \ge \Big( 1 - \frac{1}{|A|} \Big) t^2
\ge \frac{t^2}{2}
$$
provided that $|A| > 1$.

Together with \eqref{exys} this proves the existence of a
coordinate $i \in \Omega$, on which
\begin{equation}                                  \label{sxi}
  \s(x(i))  \ge  \frac{t}{2},
\end{equation}
and the claim follows from Lemma \ref{l:small deviation} applied
to the random variable $x(i)$.
\end{proof}

This lemma should be interpreted as a separation lemma for the set
$A$. It means that one can always find two nontrivial subsets of
$A$ and a coordinate in $\Omega$, on which the two subsets are
separated with a ``gap" proportional to $t$.

Based on Lemma \ref{separationlemma}, one can construct a large
separating tree in $A$. Recall that a {\em tree of subsets} of a
set $A$ is a finite collection $T$ of subsets of $A$ such that,
for every pair $B, D \in T$ either $B$ and $D$ are disjoint or one
of them contains the other. We call $D$ a {\em son} of $B$ if $D$
is a maximal (with respect to inclusion) proper subset of $B$ that
belongs to $T$. An element of $T$ with no sons is called a {\em
leaf}.

\begin{definition}
  Let $A$ be a class of functions on $\Omega$ and $t > 0$.
  A {\em $t$-separating tree $T$ of $A$} is a tree of subsets of $A$
  such that every element $B \in T$ which is not a leaf has
  exactly two sons $B_+$ and $B_{-}$ and, for some
  coordinate $i \in \Omega$,
  $$
  f(i) > g(i) + t
  \ \ \ \text{for all $f \in B_+$,  $g \in B_-$.}
  $$
\end{definition}

\begin{proposition}                      \label{thm:separatingtree}
Let $A$ be a finite class of functions on a probability space
$(\Omega,\mu)$. If $A$ is $t$-separated with respect to the
$L_2(\mu)$ norm, then there exists a $\frac{1}{6} t$-separating
tree of $A$ with at least $|A|^{1/2}$ leaves.
\end{proposition}

\begin{proof}
By Lemma \ref{separationlemma}, any finite class $A$ which is
$t$-separated with respect to the $L_2(\mu)$ norm has two subsets
$A_+$ and $A_{-}$ and a coordinate $i \in \Omega$ for which $f(i)
> g(i) + \frac{1}{6} t$ for every $f \in A_+$ and $g \in A_{-}$.
Moreover, there exists some number $0<\beta \leq 1/2$ such that
\begin{equation*}
|A_+|\geq (1-\beta)|A| \ \ \ \text{and} \ \ \ |A_{-}| \geq
\frac{\beta}{2}, \ \ \ \text{or vice versa}.
\end{equation*}
Thus, $A_+$ and $A_{-}$ are sons of $A$ which are both large and
well separated on the coordinate $i$.

The conclusion of the proposition will now follow by induction on
the cardinality of $A$. The proposition clearly holds for $|A|=2$.
Assume it holds for every $t$-separated class of cardinality
bounded by $N$, and let $A$ be a $t$-separated class of
cardinality $N+1$. Let $A_+$ and $A_{-}$ be the sons of $A$ as
above; since $\b > 0$, we have $|A_+|,|A_{-}| \leq N$. Moreover,
if $A_+$ has a $\frac{1}{6} t$-separating tree with $N_+$ leaves
and $A_{-}$ has a $\frac{1}{6} t$-separating tree with $N_{-}$
leaves then, by joining these trees, $A$ has a $\frac{1}{6}
t$-separating tree with $N_+ + N_{-}$ leaves, the number bounded
below by $|A_+|^{1/2}+|A_-|^{1/2}$ by the induction hypothesis.
Since $\b \le 1/2$,
\begin{align*}
|A_+|^{\frac{1}{2}}+|A_{-}|^{\frac{1}{2}}
  &\geq  \bigl((1-\beta)|A|\bigr)^{\frac{1}{2}}
         +\bigl(\frac{\beta}{2}|A|\bigr)^{\frac{1}{2}}\\
  &=  \Bigl[(1-\beta)^{\frac{1}{2}}
     +\bigl(\frac{\beta}{2}\bigr)^{\frac{1}{2}}\Bigr]
     |A|^{\frac{1}{2}}
  \geq |A|^{\frac{1}{2}}
\end{align*}
as claimed.
\end{proof}

The exponent $1/2$ has no special meaning in Proposition
\ref{thm:separatingtree}. It can be improved to any number smaller
that $1$ at the cost of reducing the constant $\frac{1}{6}$.

\subsection*{Counting shattered sets}
%......................................................................

As explained in the introduction, our aim is to construct a large
set shattered by a given class. We will first try to do this for
classes of integer-valued functions.

Let $A$ be a class of integer-valued functions on a set $\Omega$.
We say that a couple $(\s, h)$ is a {\em center} if $\s$ is a
finite subset of $\Omega$ and $h$ is an integer-valued function on
$\s$. We call the cardinality of $\s$ the dimension of the center.
For convenience, we introduce (the only) $0$-dimensional center
$(\emptyset, \emptyset)$, which is the {\em trivial center}.

\begin{definition}                  \label{shatteredcenter}
The {\em set $A$ shatters a center $(\s, h)$} if the following
holds:
\begin{itemize}
  \item{either $(\s, h)$ is trivial and $A$ is nonempty,}
  \item{or, otherwise, for every choice of signs
        $\theta \in \{-1,1\}^\s$
        there exists a function $f \in A$ such that for $i \in \s$
        \begin{equation}
\label{discrete shatter}
        \begin{cases}
          f(i)  > h(i)    & \text{when $\theta(i) = 1$},\\
          f(i)  < h(i)    & \text{when $\theta(i) = -1$}.
        \end{cases}
        \end{equation}
        }
\end{itemize}
\end{definition}
It is crucial that both inequalities in \eqref{discrete shatter}
are strict: they ensure that whenever a $d$-dimensional center is
shattered by $A$, one has $\vc(A, 2) \ge d$. In fact, it is
evident that $\vc(A,2)$ is the maximal dimension of a center
shattered by $A$.

\begin{proposition}                        \label{centers vs leaves}
  The number of centers shattered by $A$ is
  at least the number of leaves in any $1$-separating tree of $A$.
\end{proposition}

\proof Given a class $B$ of integer-valued functions, denote by
$s(B)$ the number of centers shattered by $B$. It is enough to
prove that if $B_+$ and $B_-$ are the sons of an element $B$ of a
$1$-separating tree in $A$ then
\begin{equation}                                   \label{HHH}
  s(B)  \ge  s(B_+) + s(B_-).
\end{equation}
By the definition of the $1$-separating tree, there is a
coordinate $i_0 \in \Omega$, such that $f(i_0) > g(i_0) + 1$ for
all $f \in B_+$ and $g \in B_-$. Since the functions are
integer-valued, there exists an integer $t$ such that
$$
f(i_0) > t \ \ \text{for $f \in B_+$} \ \ \ \text{and} \ \ \
g(i_0) < t \ \ \text{for $g \in B_-$}.
$$

If a center $(\s, h)$ is shattered either by $B_+$ or by $B_-$, it
is also shattered by $B$. Next, assume that $(\s, h)$ is shattered
by both $B_+$ and $B_-$. Note that in this case $i_0 \not\in \s$.
Indeed, if the converse holds then $\s$ contains $i_0$ and hence
is nonempty. Thus the center $(x,\s)$ is nontrivial and there
exist $f \in B_+$ and $g \in B_-$ such that $t < f(i_0) < h(i_0)$
(by \eqref{discrete shatter} with $\theta(i_0) = -1$) and $t >
g(i_0) > h(i_0)$ (by \eqref{discrete shatter} with $\theta(i_0) =
1$), which is impossible. Consider the center $(\s', h') = (\s
\cup \{i_0\}, h \oplus t)$, where $h \oplus t$ is the extension of
the function $h$ onto the set $\s \cup \{i_0\}$ defined by $(h
\oplus t) (i_0) = t$.

Observe that $(\s', h')$ is shattered by $B$. Indeed, since $B_+$
shatters $(\s, h)$, then for every $\theta \in \{-1,1\}^\s \times
\{1\}^{\{i_0\}}$ there exists a function $f \in B_+$ such that
\eqref{discrete shatter} holds for $i \in \s$. Also, since $f \in
B_+$, then automatically $f(i_0) > t = h'(i_0)$. Similarly, for
every $\theta \in \{-1,1\}^\s \times \{-1\}^{\{i_0\}}$, there
exists a function $f \in B_-$ such that \eqref{discrete shatter}
holds for $i \in \s$ and automatically $f(i_0) < t = h'(i_0)$.

Clearly, $(\s', h')$ is shattered by neither $B_+$ nor by $B_-$,
because $f(i_0) > t = h'(i_0)$ for all $f \in B_+$, so
\eqref{discrete shatter} fails if $\theta(i_0) = -1$; a similar
argument holds for $B_-$.

Summarizing, $(\s, h) \to  (\s', h')$ is an injective mapping from
the set of centers shattered by both $B_+$ and $B_-$ into the set
of centers shattered by $B$ but not by $B_+$ or $B_-$, which
proves our claim.
\endproof

Combining  Propositions \ref{thm:separatingtree} and \ref{centers
vs leaves}, one bounds from below the number of shattered centers.

\begin{corollary}                               \label{c: lowerbound}
  Let $A$ be a finite class of integer-valued functions on
  a probability space $(\Omega,\mu)$.
  If $A$ is $6$-separated with respect to the $L_2(\mu)$ norm
  then it shatters at least $|A|^{1/2}$ centers.
\end{corollary}

To show that there exists a large dimensional center shattered by
$A$, one must assume that the class $A$ is bounded in some sense,
otherwise one could have infinitely many low dimensional centers
shattered by the class. A natural assumption is the uniform
boundedness of $A$, under which we conclude a preliminary version
of Theorem \ref{main}.

\begin{proposition}                             \label{p: upperbound}
  Let $(\Omega,\mu)$ be a probability space,
  where $\Omega$ is a finite set of cardinality $n$.
  Assume that $A$ is a class of functions on
  $\Omega$ into $\{0, 1, \ldots, p \}$,
  which is $6$-separated in $L_2(\mu)$.
  Set $d$ to be the maximal dimension of a center
  shattered by $A$.
  Then
  \begin{equation}                              \label{A discrete}
    |A|  \le  \Big( \frac{p n}{d} \Big)^{C d},
  \end{equation}
  where $C$ is an absolute constant.
  In particular, the same assertion holds for $d=\vc(A,2)$.
\end{proposition}

\proof By Corollary \ref{c: lowerbound}, $A$ shatters at least
$|A|^{1/2}$ centers. On the other hand, the total number of
centers whose dimension is at most $d$ that a class of $\{0, 1,
\ldots, p \}$-valued functions on $\Omega$ can shatter is bounded
by $\sum_{k=0}^d \binom{n}{k} p^k$. Indeed, for every $k$ there
exist at most $\binom{n}{k}$ subsets $\s \subset \Omega$ of
cardinality $k$ and, for each $\s$ with $|\s| = k$ there are at
most $p^k$ level functions $h$ for which the center $(\s, h)$ can
be shattered by such a class. Therefore $|A|^{1/2} \le
\sum_{k=0}^d \binom{n}{k} p^k$ (otherwise there would exist a
center of dimension larger than $d$ shattered by $A$,
contradicting the maximality of $d$). The proof is completed by
approximating the binomial coefficients using Stirling's formula.
\endproof

Actually, the ratio $n/d$ can be eliminated from \eqref{A
discrete} (perhaps at the cost of increasing the separation
parameter $6$). To this end, one needs to reduce the size of
$\Omega$ without changing the assumption that the class is ``well
separated".  This is achieved by the following probabilistic
extraction principle.

\begin{lemma}                                      \label{extraction}
  There is a positive absolute constant $c$ such that the
  following holds.
  Let $\Omega$ be a finite set with the uniform probability measure
  $\mu$ on it.
  Let $A$ be a class of functions bounded by $1$, defined
  on $\Omega$.
  Assume that for some $0 < t < 1$
  $$
  \text{$A$ is $t$-separated with respect to the $L_2(\mu)$ norm.}
  $$
  If $|A| \leq  \frac{1}{2} \exp(c t^4 k)$
  for some positive number $k$,
  there exists a subset $\s \subset \Omega$
  of cardinality at most $k$ such that
  $$
  \text{$A$ is $\frac{t}{2}$-separated with respect to the
$L_2(\mu_\s)$
    norm,}
  $$
  where $\mu_\s$ is the uniform probability measure on $\s$.
\end{lemma}

As the reader guesses, the set $\s$ will be chosen randomly in
$\Omega$. We will estimate probabilities using a version of
Bernstein's inequality (see e.g. \cite{VW}, or \cite{LT} 6.3 for
stronger inequalities).

\begin{lemma}[Bernstein's inequality]
\label{thm:bernstein} Let $X_1, \ldots, X_n$ be independent random
variables with zero mean. Then, for every $u>0$,
\begin{equation*}
  \P \Big\{ \big| \sum_{i=1}^n X_i \big| > u \Big\}
  \leq 2 \exp \Big(-\frac{u^2}{2(b^2 + a u/3)} \Big),
\end{equation*}
where $a=\sup_i \|X_i\|_\infty$ and $b^2=\sum_{i=1}^n \E |X_i|^2$.
\end{lemma}

\noindent {\bf Proof of Lemma \ref{extraction}. } For the sake of
simplicity we identify $\Omega$ with $\{1,2, \ldots, n\}$. The
difference set $S = \{f - g |\; f \not = g, \ f,g \in A\}$ has
cardinality $|S| \leq |A|^2$. For each $x \in S$ we have $|x(i)|
\leq 2$ for all $i \in \{1,...,n\}$ and $\sum_{i=1}^n |x(i)|^2
\geq t^2 n$. Fix an integer $k$ satisfying the assumptions of the
lemma and let $\delta_1, \ldots, \d_n$ be independent
$\{0,1\}$-valued random variables with $\E \delta_i = \frac{k}{2n}
=: \delta$. Then for every $z \in S$
\begin{align*}
\P \Big\{\sum_{i=1}^n \delta_i |x(i)|^2 \leq \frac{t^2\delta n}{2}
\Big\}
  & \leq \P \Big\{ \Big| \sum_{i=1}^n \delta_i |x(i)|^2
                - \delta \sum_{i=1}^n |x(i)|^2 \Big |
                 > \frac{t^2 \delta n}{2} \Big\} \\
  &=  \P \Big\{ \Big|\sum_{i=1}^n (\delta_i -\delta)|x(i)|^2 \Big|
              > \frac{t^2 \d n}{2} \Big\} \\
  &\le 2 \exp \Big(-\frac{c t^4 \d n}{1+t^2} \Big)
   \le  2 \exp (-c t^4 k),
\end{align*}
where the last line follows from Bernstein's inequality for
$a=\sup_i \|X_i\| \leq 2$ and
$$
b^2 = \sum_{i=1}^n \E |X_i|^2
    = \sum_{i=1}^n |x(i)|^4 \; \E(\delta_i - \delta)^2 \leq 16 \d n.
$$
Therefore, by the assumption on $k$
$$
\P \Big\{ \exists x \in S :
  \Big( \frac{1}{k} \sum_{i = 1}^n \d_i |x(i)|^2 \Big)^{1/2}
      \leq \frac{t}{2} \Big\}
\le |S| \cdot 2 \exp(- c t^4 k) < 1/2.
$$
Moreover, if $\s$ is the random set $\{i \,|\; \delta_i =1\}$ then
by Chebyshev's inequality,
$$
\P \{ |\s| > k \} = \P \Big\{ \sum_{i=1}^n \d_i > k \Big\} \le
1/2,
$$
which implies that
$$
\P \big \{ \exists x \in S : \|x\|_{L_2(\mu_\s)} \le \frac{t}{2}
\big\} < 1.
$$
This translates into the fact that with positive probability the
class $A$ is $\frac{t}{2}$-separated with respect to the
$L_2(\mu_\s)$ norm.
\endproof

\qquad

\noindent {\bf Proof of Theorem \ref{main}. } One may clearly
assume that $|A| > 1$ and that the functions in $A$ are defined on
a finite domain $\Omega$, so that the probability measure $\mu$ on
$\Omega$ is supported on a finite number of atoms. Next, by
splitting these atoms (by replacing an atom $\w$ by, say, two
atoms $\w_1$ and $\w_2$, each carrying measure $\frac{1}{2}
\mu(\w)$ and by defining $f(\w_1) = f(\w_2) = f(\w)$ for $f \in
A$), one can make the measure $\mu$ almost uniform without
changing neither the covering numbers nor the shattering dimension
of $A$. Therefore, assume that the domain $\Omega$ is $\{1, 2,
\ldots, n\}$ for some integer $n$, and that $\mu$ is the uniform
measure on $\Omega$.

Fix $0 < t \leq 1/2$ and let $A$ be a $2 t$-separated in the
$L_2(\mu)$ norm. By Lemma \ref{extraction}, there is a set of
coordinates $s \subset \{1,...,n\}$ of size $|\s| \leq
\frac{C\log|A|}{t^4}$ such that $A$ is $t$-separated in
$L_2(\mu_\s)$, where $\mu_\s$ is the uniform probability measure
on $\s$.

Let $p = \lfloor 7/t \rfloor$, define $\tilde{A} \subset
\{0,1,...,p\}^\s$ by
$$
\tilde{A}=\Bigl\{ \Bigl( \Bigl \lfloor\frac{7f(i)}{t} \Big \rfloor
\Bigr)_{i \in \s} \, | \; f \in A \Bigr\},
$$
and observe that $\tilde{A}$ is $6$-separated in $L_2(\mu_\s)$. By
Proposition \ref{p: upperbound},
   $$
    |A|=|\tilde{A}|  \le  \Big( \frac{p |\sigma|}{d} \Big)^{C d}
   $$
where $d=\vc(\tilde{A},2)$, implying that
   $$
    |A| \leq \Big(\frac{C\log |A|}{dt^5}\Big)^{Cd}.
   $$
By a straightforward computation,
   $$
   |A| \leq \Bigl(\frac{1}{t}\Bigr)^{Cd},
   $$
and our claim follows from the fact that $\vc(\tilde{A},2) \leq
\vc(A, t/7)$.
\endproof

\qquad

\remark Theorem \ref{main} also holds for the $L_p(\mu)$ covering
numbers for all $0 < p < \infty$, with constants $K$ and $c$
depending only on $p$. The only minor modification of the proof is
in Lemma \ref{l:xx'}, where the equations would be replaced by
appropriate inequalities.

\qquad

\section{Applications: Gaussian Processes and Convexity}
\label{s:convexity}

The first application is a bound on the expectation of the
supremum of a Gaussian processes indexed by a set $A$. Such a
bound is provided by Dudley's integral in terms of the $L_2$
entropy of $A$; the entropy, in turn, can be majorized through
Theorem \ref{main} by the shattering dimension of $A$. The
resulting integral inequality improves the main result of
M.~Talagrand in \cite{T 92}.

If $A$ be a class of functions on the finite set $I$, then a
natural Gaussian process $(X_a)_{a \in A}$ indexed by elements of
$A$ is
$$
X_a  =  \sum_{i \in I} g_i \, a(i)
$$
where $g_i$ are independent standard Gaussian random variables.

\begin{theorem} \label{thm:talagrand}
  Let $A$ be a class of functions bounded by $1$,
  defined on a finite set $I$ of cardinality $n$.
  Then $E = \E \sup_{a \in A} X_a$
  is bounded as
  $$
  E \le  K \sqrt{n}
       \int_{cE/n}^{1} \sqrt{\vc(A,t) \cdot \log (2/t)}\; dt,
  $$
  where $K$ and $c$ are absolute positive constants.
\end{theorem}

The nonzero lower limit in the integral will play an important
role in the application to Elton's Theorem.

The first step in the proof is to view $A$ as a subset of $\R^n$.
Dudley's integral inequality can be stated as
$$
E  \le  K \int_0^\infty \sqrt{\log N(A, t D_n)} \; dt,
$$
where $D_n$ is the unit Euclidean ball in $\R^n$, see \cite{Pi}
Theorem 5.6. The lower limit in this integral can be improved by a
standard argument. This fact was first noticed by A. Pajor.

\begin{lemma}                       \label{dudley}
  Let $A$ be a subset of $\R^n$.
  Then $E = \E \sup_{a \in A} X_a$
  is bounded as
  $$
  E  \le  K \int_{cE/\sqrt{n}}^\infty \sqrt{\log N(A, t D_n)} \; dt,
  $$
  where $K$ is an absolute constant.
\end{lemma}

\proof Fix positive absolute constants $c_1, c_2$ whose values
will be specified later. There exists a subset $\NN$ of $A$, which
is a $(\frac{c_1 E}{\sqrt{n}})$-net of $A$ with respect to the
Euclidean norm and has cardinality $|\NN| \le N(A, \frac{c_1 E}{2
\sqrt{n}} D_n)$. Then $A  \subset  \NN + \frac{c_1 E}{2 \sqrt{n}}
D_n$, and one can write
\begin{equation}                                       \label{E}
  E  =  \E \sup_{a \in A} X_a
    \le \E \max_{a \in \NN} X_a
        + \E \sup_{a \in \frac{c_1 E}{\sqrt{n}} D_n} X_a.
\end{equation}
The first summand is estimated by Dudley's integral as
\begin{equation}                                       \label{first NN}
  \E \max_{a \in \NN} X_a
  \le  K \int_0^\infty \sqrt{\log N(\NN, t D_n)} \; dt.
\end{equation}
On the interval $(0, \frac{c_2 E}{\sqrt{n}})$,
\begin{align*}
K \int_0^{ \frac{c_2 E}{\sqrt{n}} } \sqrt{\log N(\NN, t D_n)} \;
dt
  &\le  K \frac{c_2 E}{\sqrt{n}} \cdot \sqrt{\log|\NN|}  \\
  &\le  K \frac{c_2 E}{\sqrt{n}} \cdot
     \sqrt{\log N(A, {\textstyle \frac{c_1 E}{2 \sqrt{n}} } D_n)}.
\end{align*}
The latter can be estimated using Sudakov's inequality
\cite{D,Pi}, which states that $\e \sqrt{\log(N, \e D_n)} \le K
\,\E \sup_{a \in A} X_a$ for all $\e
> 0$. Indeed,
$$
K \frac{c_2 E}{\sqrt{n}} \cdot
     \sqrt{\log N(A, {\textstyle \frac{c_1 E}{2 \sqrt{n}} } D_n)}
     \leq
K_1 (2c_2 / c_1) \,\E \sup_{a \in A} X_a
  =  K_1 (2c_2 / c_1) E
   \le \frac{1}{4} E,
$$
if we select $c_2$ as $c_2 = c_1 / 8 K_1$. Combining this with
\eqref{first NN} implies that
\begin{equation}                              \label{first sum}
  \E \max_{x \in \NN} X_a
  \le \frac{1}{4} E + K \int_{ \frac{c_2 E}{\sqrt{n}} }^\infty
      \sqrt{\log N(A, t D_n)} \; dt
\end{equation}
because $\NN$ is a subset of $A$.

To bound the second summand in \eqref{E}, we apply the
Cauchy-Schwarz inequality to obtain that for any $t>0$,
$$
\E \sup_{a \in tD_n} X_a \le t \cdot \E \Big( \sum_{i \in I} g_i^2
\Big)^{1/2} \le t \sqrt{n}.
$$

In particular, if $c_1 < 1/4$ then
$$
\E \sup_{a \in \frac{c_1E}{\sqrt{n}}D_n} X_a \leq c_1E \leq
\frac{1}{4}E.
$$
This, \eqref{E} and \eqref{first sum} imply that
$$
E \le K_2 \int_{ \frac{c_2 E}{\sqrt{n}} }^\infty
      \sqrt{\log N(A, t D_n)} \; dt,
$$
where $K_2$ is an absolute constant.
\endproof

\qquad

\noindent{\bf Proof of Theorem \ref{thm:talagrand}. } By Lemma
\ref{dudley},
$$
E  \le  K \int_{cE/\sqrt{n}}^\infty \sqrt{\log N(A, t D_n)} \; dt.
$$
Since $A  \subset [-1,1]^n \subset \sqrt{n} D_n$, the integrand
vanishes for $t \ge \sqrt{n}$. Hence, by Theorem \ref{main}
\begin{align*}
E &\le  K \int_{cE/\sqrt{n}}^{\sqrt{n}}
          \sqrt{\log N(A, t D_n)} \; dt \\
  &=    K \sqrt{n} \int_{cE/n}^1
          \sqrt{\log N(A, t \sqrt{n} D_n)} \; dt \\
  &\le  K_1\sqrt{n} \int_{cE/n}^1
          \sqrt{\vc(A, c_1 t) \cdot \log (2/t)} \; dt.
\end{align*}
The absolute constant $0 < c_1 < 1/2$ can be made $1$ by a further
change of variable.
\endproof

\qquad

The main consequence of Theorem \ref{thm:talagrand} is Elton's
Theorem with the optimal dependence on $\d$.

\begin{theorem}                             \label{thm:elton}
  There is an absolute constant $c$ for which the following
  holds.
  Let $x_1, \ldots, x_n$ be vectors in the unit
  ball of a Banach space.
  Assume that
  $$
  \E  \Big\| \sum_{i = 1}^n g_i x_i \Big\|  \ge  \d n
  \ \ \ \text{for some number $\d > 0$}.
  $$
  Then there exist numbers $s, t \in (c \d, 1)$,
  and a subset $\s \subset \{ 1, \ldots, n \}$
  of cardinality $|\s|  \ge  s^2 n$,
  such that
  \begin{equation}                      \label{lower l1}
    \Big\| \sum_{i \in \s} a_i x_i \Big\|
    \ge   t \sum_{i \in \s} |a_i|
    \ \ \ \text{for all real numbers $(a_i)$}.
  \end{equation}
  In addition, the numbers $s$ and $t$ satisfy the inequality
  $s \cdot t \log^{1.6} (2 / t)  \ge  c\d$.
\end{theorem}

Before the proof, recall the interpretation of the shattering
dimension of convex bodies. If a set $B \subset \R^n$ is convex
and symmetric then $\vc(B,t)$ is the maximal cardinality of a
subset $\s$ of $\{1, \ldots, n\}$ such that $P_\s(B) \supset
[-\frac{t}{2}, \frac{t}{2}]^\s$. Indeed, every convex symmetric
set in $\R^n$ can be viewed as a class of functions on
$\{1,...,n\}$. If $\s$ is $t$-shattered with a level function $h$
then for every $\s' \subset \s$ there is some $f_{\s'}$ such that
$f_{\s'}(i) \geq h(i)+t$ if $i \in \s'$ and $f_{\s'} \leq h$ on
$\s \backslash \s'$. By selecting for every such $\s'$ the
function $(f_{\s'}-f_{\s \backslash \s'})/2$ and since the class
is convex and symmetric, it follows that $P_\s(B) \supset
[-\frac{t}{2},\frac{t}{2}]^\s$, as claimed.

Taking the polars, this inclusion can be written as $\frac{t}{2}
(B^\circ \cap \R^\s) \subset B_1^n$, where $B_1^n$ is the unit
ball of $\ell_1^n$. Denoting by $\|\cdot\|_{B^\circ}$ the
Minkowski functional (the norm) induced by the body $B^\circ$, one
can rewrite this inclusion as the inequality
$$
\Big\| \sum_{i \in \s} a_i e_i \Big\|_{B^\circ}
  \ge   \frac{t}{2} \sum_{i \in \s} |a_i|
  \ \ \ \text{for all real numbers $(a_i)$},
$$
where $(e_i)$ is the standard basis of $\R^n$. Therefore, to prove
Theorem \ref{thm:elton}, one needs to bound below the shattering
dimension of the dual ball of a given Banach space.

\qquad

\noindent {\bf Proof of Theorem \ref{thm:elton}.} By a
perturbation argument, one may assume that the vectors $(x_i)_{i
\le n}$ are linearly independent. Hence, using an appropriate
linear transformation one can assume that $X = (\R^n, \|\cdot\|)$
and that $(x_i)_{i \le n}$ are the unit coordinate vectors
$(e_i)_{i \le n}$ in $\R^n$. Let $B=(B_X)^\circ$ and note that the
assumption $\|e_i\|_X \leq 1$ implies that $B \subset [-1,1]^n$.

Set
$$
E = \E \Big\|\sum_{i=1}^n g_i x_i \Big\|_X
  = \E \sup_{b \in B} \sum_{i=1}^n g_i \,b(i).
$$
By Theorem \ref{thm:talagrand},
\begin{equation*}
\d n \leq E \leq K\sqrt{n} \int_{c \delta}^{1}
               \sqrt{\vc(B,t) \cdot \log(2/t)} \; dt.
\end{equation*}
Consider the function
$$
h(t)  =  \frac{c_0}{t \log^{1.1} (2 / t)}
$$
where the absolute constant $c_0 > 0$ is chosen so that
$\int_0^{1} h(t) \; dt  =  1$. It follows that there exits some
$c\d  \le  t \le  1$ such that
$$
\sqrt{\vc(B, t) / n \cdot \log(2 / t)} \ge  \d h(t).
$$
Hence
$$
\vc(B, t)  \ge  \frac{c_0 \d^2}{t^2 \log^{3.2} (2 / t)} n.
$$
Therefore, letting $s^2 = \vc(B, t) / n$, it follows that $s \cdot
t \log^{1.6} (2 / t)  \ge  \sqrt{c_0} \d$ as required, and by the
discussion preceding the proof there exists a subset $\s$ of $\{1,
\ldots, n\}$ of cardinality $|\s|  \ge  s^2 n$ such that
\eqref{lower l1} holds with $t/2$ instead of $t$. The only thing
remaining is to check that $s \gtrsim \d$. Indeed, $s \ge
\frac{\sqrt{c_0} \d}{t \log^{1.6} (2 / t)} \ge c_1 \d$, because $t
\le 1$.
\endproof

\noindent {\bf Remarks. } 1. As the proof shows, the exponent
$1.6$ can be reduced to any number larger than $3/2$.

2. The relation between $s$ and $t$ in Theorem \ref{thm:elton} is
optimal up to a logarithmic factor for all $0 < \d < 1$. This is
seen from by the following example, shown to us by Mark Rudelson.
For $0 < \d < 1 / \sqrt{n}$, the constant vectors $x_i = \d
\sqrt{n} \cdot e_1$ in $X = \R$ show that $s t$ in Theorem
\ref{thm:elton} can not exceed $\d$. For $1 / \sqrt{n}  \le \d \le
1$, we consider the body $D = \conv( B_1^n \cup \frac{1}{\d
\sqrt{n}} D_n )$ and let $X = (\R^n, \|\cdot\|_D)$ and $x_i =
e_i$, $i = 1, \ldots, n$. Clearly, $\E \| \sum g_i x_i \|_X  \ge
\E \| \sum \e_i e_i \|_D = \d n$. Let $0 < s, t < 1$ be so that
\eqref{lower l1} holds for some subset $\s \subset \{ 1, \ldots, n
\}$ of cardinality $|\s|  \ge  s^2 n$. This means that $\|x\|_D
\ge  t \|x\|_1$ for all $x \in \R^\s$. Dualizing, $\frac{t}{\d
\sqrt{n}} \|x\|_2  \le  t \|x\|_{D^\circ}  \le \|x\|_\infty$ for
all $x \in \R^\s$. Testing this inequality for $x = \sum_{i \in
\s} e_i$, it is evident that $\frac{t}{\d \sqrt{n}} \sqrt{|\s|}
\le 1$ and thus $s t  \le  \d$.

\qquad

We end this article with an application to empirical processes. A
key question is when a class of functions satisfies the central
limit theorem uniformly in some sense. Such classes of functions
are called {\it uniform Donsker classes}. We will not define these
classes formally but rather refer the reader to \cite{D,VW} for an
introduction on the subject. It turns out that the uniform Donsker
property is related to uniform estimates on covering numbers via
the Koltchinskii-Pollard entropy integral.

\begin{theorem} \cite{D}
\label{dud} Let $F$ be a class of functions bounded by $1$. If
\begin{equation*}
 \int_0^\infty  \sup_{n}\sup_{\mu_n}  \sqrt{ \log N
\bigl(F,L_2(\mu_n),\e\bigr) } \;d\eps < \infty,
\end{equation*}
then $F$ is a uniform Donsker class.
\end{theorem}
Having this condition in mind, it is natural to try to seek
entropy estimates which are ``dimension free", that is, do not
depend on the size of the sample. In the $\{0,1\}$-valued case,
such bounds where first obtained by Dudley who proved Theorem
\ref{main} for these classes (see \cite{LT} Theorem 14.13) which
implied through Theorem \ref{dud} that every VC class is a uniform
Donsker class.

Theorem \ref{main} solves the general case: the following
corollary extends Dudley's result on the uniform Donsker property
from $\{0,1\}$ classes to classes of real valued functions.

\begin{corollary}
  Let $F$ be a class of functions bounded by $1$
  and assume that the integral
  $$
  \int_0^1 \sqrt{ \vc(F,t)\log \frac{2}{t}} \; dt
  $$
  converges.
  Then $F$ is a uniform Donsker class.
\end{corollary}

In particular this shows that if $\vc(F,t)$ is ``slightly better"
than $1/t^2$, then $F$ is a uniform Donsker class.

This result has an advantage over Theorem \ref{dud} because in
many cases it is easier to compute the shattering dimension of the
class rather than its entropy (see, e.g. \cite{AB}).

\end{document}